# Study of memory effect in an EOQ model for completely backlogged demand during shortage

Rituparna Pakhira[a], Uttam Ghosh[b], Susmita Sarkar[c], Vishnu Narayan Mishra[d]

[a,b,c]Department of Applied Mathematics, University of Calcutta, Kolkata, West Bengal, India

[d]Department of Applied Mathematics, Indira Gandhi National Tribal University, Lalpur, Amarkantak 484 887, Anuppur, Madhya

[a]email:rituparna.pakhira@gmail.com
[b]email:uttam_math@yahoo.com
[c]email:susmita62@yahoo.co.in
[d]email:vishnunarayanmishra@gmail.com

**Abstract**

The most commonly developed inventory models are the classical economic order quantity model, is governed by the integer order differential equations. We want to come out from the traditional thought *i.e.* classical order inventory model where the memory phenomena are absent. Here, we want to incorporate the memory effect that is based on the fact economic agents remember the history of changes of exogenous and endogenous variables. In this paper, we have proposed and solved a fractional order EOQ model with constant demand rate where the demand is fully backlogged during shortage time. Finally, a numerical example has been illustrated for this model to show the memory dependency of the system. The numerical example clears that for the considered system the profit is maximum in long memory affected system compared to the low memory affected or memory less system.

1. **Introduction**

Recently, fractional order integration and differentiation have been applied to the real world problem for its memory property. But in this system growth of any processes is slower compared to the ordinary differential system. Fractional order integration and differentiation are generalizations of integer-order integration and differentiation with including $n^{th}$ derivatives and fractional $n$-fold integrals [8-9]. Usually, authors use R-L integration and Caputo fractional order derivative to develop different problems [1-8]. Recently Baleanu and his collaborators developed and used the new non-local fractional derivative with non-singular kernel in different physical problems [26-28]. The fractional order derivative has different type of significance in different applications. In physical problem it is used as the roughness parameter of the surface [21], in biological and financial system it is the indicator of memory [1-8] etc.

The memory means the dependence of the process not only on the current state of the processes but also on the past history of the process. Like the application in physics and biology it is recently using in Economic analysis [2-4]. However, there are some areas of operation research where the application of fractional calculus has been started to apply in the last few years. Our main interest is to include the application fractional calculus in the inventory models.

It is well known that rate of change of integer orders of a function at any particular point is determined by the property of the function in the infinitely small neighborhood of the considered point. Hence, the integer order rate of change of any function/system is assumed as the instant rate of change of the marginal output, when the input level changes. Therefore, dynamic memory effect is not present in classical calculus. Thus the classical calculus is not able to include the previous state of the system [6-7]. But in the fractional derivative, the rate of change is affected by all points of the considered interval, so it is able to incorporate the previous /memory effects of any system. Here, the order of fractional derivative will be treated as an index of memory. So, the fractional order system can remove amnesia from any system [22].

Its application is suitable to formulate and analyze the real-life phenomena which has memory effects. Inclusion of memory effect using fractional calculus in susceptible-infected-recovered (SIR) epidemic model was done by M.Saedian et al. [1] using the memory kernel function. They have proposed in their paper [1] that the real epidemic process is clearly sustained by a non-Markovian dynamics memory effect.

Tarasova et al.[2-5] have developed many research articles using the concept of the memory effect of fractional order derivative and integration. In the paper [4] Tarasova et al. gave an idea of economic interpretation of fractional derivatives using the Caputo derivative. Tarasova et al. in [5] discussed elasticity for the economic process with memory using fractional differential calculus. Here, they defined that generalization of point price elasticity of demand to the case of the processes with memory. In these generalizations, they take into account dependence of demand not only from the current price(price at the current time) but also changes of price for some time interval.

Pakhira et al.[6-7, 23-24] developed some memory dependent Economic order quantity models including the fractional order rate of change of the inventory level and fractional order effect of different costs. The inclusion of the memory effect in the inventory model is necessary to handle practical business policy.

In this paper, our aim is to develop a memory dependent EOQ model where demand is completely backlogged during shortage [25]. Here, Caputo fractional order derivative and Riemann-Liouville fractional order integration have been used to develop the EOQ model where fully backlogged is permitted during shortage time. To formulate the fractional order differential equation model we have used the memory kernel concept as develop in [ 1]. Here different costs are established using fractional order integration. Finally the effects of inclusion of memory parameter are investigated using numerical examples.

Our analysis clears that when the integral memory index is absent but the differential memory index is present, the minimized total average cost is gradually decreasing with gradually increasing memory effect but optimal ordering interval is gradually increasing. The numerical value of the minimized total average cost with the presence of the differential memory index is low compared to the presence of the integral memory index. Analysis also clear that total order quantity needs to be adjusted for short past experience effect but not for long experience.

The rest part of the paper is organized as follows in the section-2, some review of fractional calculus has been given, in the section-3, model formulation has been discussed, in the section-3.4, fractional order inventory model formulation with memory kernel has been given, In the section -3.5 fractional order model Analysis has been presented, In the section-4, numerical examples are given and finally in the section-5, some conclusions are cited.

## 2. Review of fractional calculus

### 2.1 *Euler Gamma Function*
Euler's gamma function is one of the best tools in fractional calculus which was proposed by the Swiss mathematicians Leonhard Euler (1707-1783).The gamma function $\Gamma(x)$ is continuous extension from the factorial notation. The gamma function is denoted and defined by the formulae

$$\Gamma(x) = \int_0^\infty t^{(x-1)} e^{-t} dt \qquad x > 0 \qquad (1)$$

$\Gamma(x)$ is extended for all real and complex numbers. It has the basic properties $\Gamma(x+1) = x\Gamma(x), \Gamma\left(\frac{1}{2}\right) = \frac{\sqrt{\pi}}{2}, \Gamma\left(-\frac{7}{6}\right) = -\frac{6}{7}\Gamma\left(\frac{1}{6}\right)$ etc.

Numerically, $\lfloor x$ can be evaluated for positive integer values only $x$ but $\Gamma(x+1)$ can be evaluated for real values of $x$.

### 2.2 Riemann-Liouville fractional derivative(R-L)
Left Riemann-Liouville fractional derivative of order $\alpha$ is denoted and defined as follows

$$_aD_x^\alpha(f(x)) = \frac{1}{\Gamma(m-\alpha)} \left(\frac{d}{dx}\right)^m \int_a^x (x-\tau)^{(m-\alpha-1)} f(\tau) d\tau \qquad (2)$$

where $x > 0$

Right Riemann-Liouville fractional derivative of order $\alpha$ is defined as follows

$$_xD_b^\alpha(f(x)) = \frac{1}{\Gamma(m-\alpha)} \left(-\frac{d}{dx}\right)^m \int_x^b (x-\tau)^{(m-\alpha-1)} f(\tau) d\tau \qquad (3)$$

where $x > 0$

Riemman-Liouville fractional derivative of any constant function is not equal to zero which creates a distance between ordinary calculus and fractional calculus. This definition creates a difficulty that action of derivative of constant term is not zero.

### 2.3 Caputo fractional order derivative

Left Caputo fractional derivative [10] for the function $f(x)$ which has continuous, bounded derivatives in $[a,b]$ is denoted and defined as follows

$$_a^C D_x^\alpha (f(x)) = \frac{1}{\Gamma(m-\alpha)} \int_a^x (x-\tau)^{(m-\alpha-1)} f^m(\tau) d\tau \quad (4)$$

where $0 \leq m-1 < \alpha < m$

Right Caputo fractional derivative for the function $f(x)$ which has continuous and bounded derivatives in $[a,b]$ is defined as follows

$$_x^C D_b^\alpha (f(x)) = \frac{1}{\Gamma(m-\alpha)} \int_x^b (\tau-x)^{(m-\alpha-1)} f^m(\tau) d\tau \quad (5)$$

where $0 \leq m-1 < \alpha < m$

$_a^C D_x^\alpha (A) = o,$ where $A =$ constant.

### 2.4 Fractional Laplace transforms Method

The Laplace transform of the function $f(t)$ is defined as

$$F(s) = L(f(t)) = \int_0^\infty e^{-st} f(t) dt \quad (6)$$

where $s > 0$ and $s$ is called the transform parameter. The Laplace transformation of $n^{th}$ order derivative is defined as

$$L(f^n(t)) = s^n F(s) - \sum_{k=0}^{n-1} s^{n-k-1} f^k(0) \quad (7)$$

where $f^n(t)$ denotes $n^{th}$ derivative of the function $f$ with respect to $t$. For non-integer $\alpha$ it is defined in generalized form [8-9] as,

$$L(f^\alpha(t)) = s^\alpha F(s) - \sum_{k=0}^{n-1} s^k f^{\alpha-k-1}(0) \quad (8)$$

where, $(n-1) < \alpha \leq n$.

### 3. Model formulation

Inventory models are developed under the following assumptions and notations, are given below.

### 3.1 Assumption

In this paper, the classical and fractional order EOQ models are developed on the basis of the following assumptions.

(i) Lead time is zero.
(ii) Time horizon is infinite.
(iii) Demand rate is $D(t) = \gamma$    for $0 \leq t \leq t_1$
                      $= \gamma$    for $t_1 \leq t \leq T$.
(iv) Shortage is allowed in this model.
(v) There is permitted complete backlogging during shortage time.

### 3.2 Notations

The following notations are used to develop the model.

| $(i) D(t)$: Demand rate. | $(ii) Q$: Total order quantity. |
|---|---|
| $(iii) P$: Per unit cost. | $(iv) C_1 t^\alpha$: Inventory holding cost per unit. |
| $(v) C_3$: Ordering cost or setup cost. | $(vi) I(t)$: Stock level or inventory level |
| $(vii) T$: Ordering interval. | $(viii) HOC_{\alpha,\beta}(T)$: Inventory holding cost with fractional effect. |
| $(ix) T^*_{\alpha,\beta}$: Optimal ordering interval with fractional effect. | $(x) TOC^{av}_{\alpha,\beta}$: Total average cost during the total time interval. |
| $(xi) TOC^*_{\alpha,\beta}$: Minimized total average cost with fractional effect. | $(xii) (B,.), (\Gamma,.)$ Beta function and gamma function respectively. |
| $(xiii) SOC_{\alpha,\beta}(T)$: Shortage cost with fractional effect. | $(xiv) POC_{\alpha,\beta}(T)$: Total purchasing cost with fractional effect. |
| $(xv): C_2$: Shortage cost per unit per unit time. | |

**Table-1:** Different symbols and items for the EOQ models.

### 3.3 Classical inventory model

Here the positive inventory level $(I_1(t))$ and the negative inventory level $(I_2(t))$ both reduce due to constant demand $(\gamma)$ during the time interval $[0, t_1]$ and $[t_1, T]$ respectively. The classical or memory less inventory system with constant demand can be governed by the following system of integer order differential equations with the initial conditions:

$$\frac{d(I_1(t))}{dt} = -(\gamma) \qquad \text{for } 0 \leq t \leq t_1 \qquad (9)$$

$$\frac{d(I_2(t))}{dt} = -\gamma \qquad \text{for } t_1 \leq t \leq T \quad (10)$$

with $I_1(t_1) = I_2(t_1) = 0$

In the next section we shall establish the fractional order inventory control model using the memory kernel method [1].

### 3.4 Fractional order inventory model formulation with memory kernel

To establish the influence of memory effects, the differential equation (9-10) can be written using the kernel functions in the following form:

$$\frac{dI_1(t)}{dt} = -\int_0^t k(t-t')(\gamma) dt' \quad (11)$$

$$\frac{dI_2(t)}{dt} = -\int_0^t k(t-t') \gamma dt' \quad (12)$$

in which $k(t-t')$ is the kernel function. This type of kernel guarantees the existence of scaling features as it is often intrinsic in most natural phenomena. Thus, to generate the fractional order model we consider $k(t-t') = \frac{1}{\Gamma(\alpha-1)} (t-t')^{\alpha-2}$ where $0 < \alpha \leq 1$ and, $\Gamma(\alpha)$ denotes the gamma function.

Using the definition of fractional derivative [8-9] we can re-write the Equation (11, 12) to the form of fractional differential equations with the Caputo-type derivative in the following form

$$\frac{dI_1(t)}{dt} = -_0D_t^{-(\alpha-1)}(\gamma) \quad (13)$$

$$\frac{dI_2(t)}{dt} = -_0D_t^{-(\alpha-1)}(\gamma) \quad (14)$$

Now, applying fractional caputo derivative of order $(\alpha-1)$ on both sides of the Eq. (13, 14) and using the fact that the Caputo fractional order derivative and fractional order integral are inverse operators, the following fractional differential equations can be obtained for the model

$$_0^C D_t^\alpha (I_1(t)) = -(\gamma) \quad (15)$$

$$_0^C D_t^\alpha (I_2(t)) = -\gamma \quad (16)$$

along with boundary condition $I_1(t_1) = 0$, $I_2(t_1) = 0$

Here, $\alpha$ controls the strength of memory. When $\alpha \to 1$, memory of the system becomes weak and the small value of $\alpha$ (close to 0.1) indicates long memory of the system.

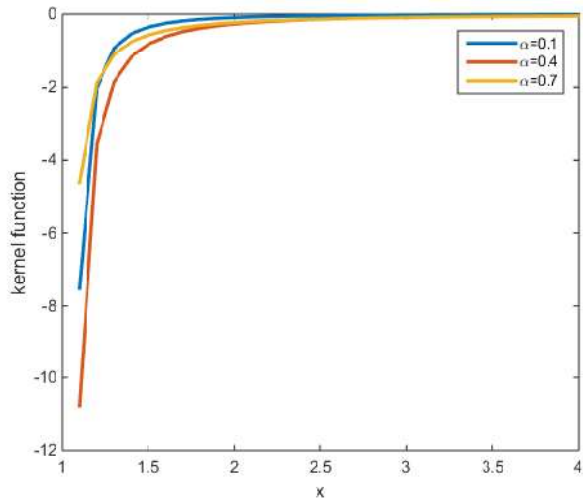

**Fig-1:** Plot of memory kernel function $k(x-\xi) = \frac{(x-\xi)^{(\alpha-2)}}{\Gamma(1-\alpha)}$, for $\xi = 3$ for different values of $\alpha$.

It is clear from the figure that the pick of the curve gradually decreases depending on $\alpha$.
Here, we define memory effect in two steps(i) long memory effect,(ii) low memory effect.

**Long and Short memory effect:**

The strength of memory is controlled by the order of fractional derivative or fractional integration.If order of fractional derivative or fractional integrationis in the range (0,0.5) then the system has long memory effect and in [0.5,1) then is short memory effect.

## 3.5 Fractional order model Analysis

Fractional order inventory model is governed by the two fractional order differential equation as follows

$$\frac{d^{\alpha}(I_1(t))}{dt^{\alpha}} = -\gamma \qquad \text{for} \qquad 0 \leq t \leq t_1 \quad (17)$$

$$\frac{d^{\alpha}(I_2(t))}{dt^{\alpha}} = -\gamma \qquad \text{for} \quad t_1 \leq t \leq T \quad (18)$$

The inventory levelcan be found integrating the system (17-18) and using the boundary condition $I_1(t_1) = I_2(t_1) = 0$ which gives

$$I_1(t) = \frac{\gamma}{\Gamma(1+\alpha)}(t_1^{\alpha} - t^{\alpha}) \quad (19)$$

$$I_2(t) = \frac{\gamma}{\Gamma(1+\alpha)}(t_1^{\alpha} - t^{\alpha}) \quad (20)$$

Since the inventory level decreases with respect to time *t*, so the maximum positive inventory level will occur at $t = 0$ which gives

$$I_1(0) = M = \left(\frac{\gamma t_1^{\alpha}}{\Gamma(\alpha+1)}\right) \quad (21)$$

Here, the maximum backorder units are

$$S = -I_2(T) = \frac{\gamma}{\Gamma(1+\alpha)}(T^{\alpha} - t_1^{\alpha}) \quad (22)$$

The order size during the total ordering interval $[0,T]$, is denoted by $Q$ and defined as

$$Q = M + S = \frac{\gamma}{\Gamma(1+\alpha)}(T^{\alpha}) \quad (23)$$

In reality, holding cost depends on time. It is not constant in the entire cycle of the system. Due to that reason, the inventory holding cost per unit is assumed as a function of time in the form $C_1 t^{\alpha}$.

The inventory holding costwith memory effect is denoted by $HOC_{\alpha,\beta}$ [6] and defined as

$$HOC_{\alpha,\beta}(T) = \frac{C_1}{\Gamma(\beta)} \int_0^{t_1} t^{\alpha}(t_1-t)^{\beta-1} I_1(t) dt = \frac{C_1 \gamma t_1^{2\alpha+\beta}(B(\alpha+1,\beta) - B(2\alpha+1,\beta))}{\Gamma(\beta)\Gamma(1+\alpha)} \quad (24)$$

($\beta$ is considered as integral memory index)

Shortage cost with fractional effect is denoted by $SC_{\alpha,\beta}$ and defined as follows

$$\begin{aligned} SOC_{\alpha,\beta}(T) &= -\frac{C_2}{\Gamma(\beta)} \int_{t_1}^{T} (T-t)^{\beta-1} I_2(t) dt \\ &= \left( \frac{C_2 \gamma T^{\alpha+\beta}}{\Gamma(\beta)\Gamma(\alpha+1)}\left(\frac{1}{(\alpha+1)} - \frac{(\beta-1)}{(\alpha+2)}\right) + \left(\frac{C_2 \gamma t_1^{\alpha+1}}{\Gamma(\beta)\Gamma(\alpha+1)} - \frac{C_2 \gamma t_1^{\alpha+1}}{(\alpha+1)\Gamma(\beta)\Gamma(\alpha+1)}\right) T^{\beta-1} \right. \\ &\left. + \left(\frac{C_2 \gamma t_1^{\alpha+2}(\beta-1)}{\Gamma(\beta)(\alpha+2)\Gamma(\alpha+1)} - \frac{C_2 \gamma t_1^{\alpha+2}(\beta-1)}{\Gamma(\beta)(2)\Gamma(\alpha+1)}\right) T^{\beta-2} + \left(\frac{C_2 \gamma(\beta-1)t_1^{\alpha}}{\Gamma(\beta)2\Gamma(\alpha+1)} - \frac{C_2 \gamma t_1^{\alpha}}{\Gamma(\beta)\Gamma(\alpha+1)}\right) T^{\beta} \right) \end{aligned} \quad (25)$$

$C_2$ is the shortage cost per unit per unit time

The purchasing cost for fractional order model is denoted by $POC_{\alpha,\beta}$ and defined as

$$POC_{\alpha,\beta}(T) = PxQ = P\left(\frac{\gamma T^{\alpha}}{\Gamma(1+\alpha)}\right) \quad (26)$$

Therefore, total average cost for the fractional order inventory model is

$$TOC_{\alpha,\beta}^{av}(T) = \frac{(HOC_{\alpha,\beta}(T) + POC_{\alpha,\beta} + SOC_{\alpha,\beta} + C_3)}{T}$$
$$= AT^{\alpha+\beta-1} + B_1 T^{\beta-1} + CT^{\beta-2} + DT^{\beta-3} + ET^{\alpha-1} + FT^{-1} \quad (27)$$

where,

$$A = \frac{C_2 \gamma}{\Gamma(\beta)\Gamma(\alpha+1)}\left(\frac{1}{(\alpha+1)} - \frac{(\beta-1)}{(\alpha+2)}\right), B_1 = \frac{C_2\gamma(\beta-1)t_1^{\alpha}}{2\Gamma(\beta)\Gamma(\alpha+1)} - \frac{C_2\gamma t_1^{\alpha}}{\Gamma(\beta)\Gamma(\alpha+1)}, C = \frac{C_2\gamma t_1^{\alpha+1}}{\Gamma(\beta)\Gamma(\alpha+1)} - \frac{C_2\gamma t_1^{\alpha+1}}{\Gamma(\beta)(\alpha+1)\Gamma(\alpha+1)}$$

$$D = \frac{C_2\gamma(\beta-1)t_1^{\alpha+2}}{(\alpha+2)\Gamma(\beta)\Gamma(\alpha+1)} - \frac{C_2\gamma(\beta-1)t_1^{\alpha+2}}{2\Gamma(\beta)\Gamma(\alpha+1)}, E = \frac{P\gamma}{\Gamma(\alpha+1)}, F = \frac{C_1\gamma(B(\alpha+1,\beta) - B(2\alpha+1,\beta))t_1^{2\alpha+\beta}}{\Gamma(\beta)\Gamma(1+\alpha)} + C_3.$$

We consider three different cases of the total average cost depending on the different values of the memory indexes $\alpha, \beta$ (i) $0 < \alpha \leq 1.0, 0 < \beta \leq 1.0$. (ii) $\beta = 1.0,\ 0 < \alpha \leq 1.0$, (iii) $\alpha = 1.0,\ 0 < \beta \leq 1.0$.

**(i) Case-1:** we consider both memory exist i.e. differential memory index and integral memory indexes both fractional i.e. $0 < \alpha \leq 1.0, 0 < \beta \leq 1.0$.

Here, the inventory model can be written as follows

$$\begin{cases} \text{Min}\, TOC_{\alpha,\beta}^{av}(T) = AT^{\alpha+\beta-1} + B_1 T^{\beta-1} + CT^{\beta-2} + DT^{\beta-3} + ET^{\alpha-1} + FT^{-1} \\ \text{Subject to}\, T \geq 0 \end{cases} \quad (28)$$

$$A = \frac{C_2\gamma}{\Gamma(\beta)\Gamma(\alpha+1)}\left(\frac{1}{(\alpha+1)} - \frac{(\beta-1)}{(\alpha+2)}\right), B_1 = \frac{C_2\gamma(\beta-1)t_1^{\alpha}}{2\Gamma(\beta)\Gamma(\alpha+1)} - \frac{C_2\gamma t_1^{\alpha}}{\Gamma(\beta)\Gamma(\alpha+1)}, C = \frac{C_2\gamma\alpha t_1^{\alpha+1}}{\Gamma(\beta)(\alpha+1)\Gamma(\alpha+1)}$$

$$D = \frac{C_2\gamma(\beta-1)t_1^{\alpha+2}}{(\alpha+2)\Gamma(\beta)\Gamma(\alpha+1)} - \frac{C_2\gamma t_1^{\alpha}}{\Gamma(\beta)\Gamma(\alpha+1)} - \frac{C_2\gamma(\beta-1)t_1^{\alpha+2}}{2\Gamma(\beta)\Gamma(\alpha+1)}, E = \frac{P\gamma}{\Gamma(\alpha+1)}, F = \frac{C_1\gamma(B(\alpha+1,\beta) - B(2\alpha+1,\beta))t_1^{2\alpha+\beta}}{\Gamma(\beta)\Gamma(1+\alpha)} + C_3.$$

**(a) Primal Geometric Programming Method**

To solve (28) analytically, the primal geometric programming method has been applied. The dual form of (28) has been introduced by the dual variable ($w$). The corresponding primal geometric programming problem has been constructed in the following:

$$\text{Max}\, d(w) = \left(\frac{A}{w_1}\right)^{w_1} \left(\frac{B_1}{w_2}\right)^{w_2} \left(\frac{C}{w_3}\right)^{w_3} \left(\frac{D}{w_4}\right)^{w_4} \left(\frac{E}{w_5}\right)^{w_5} \left(\frac{F}{w_6}\right)^{w_6} \quad (29)$$

Normalized condition is as

$$w_1 + w_2 + w_3 + w_4 + w_5 + w_6 = 1 \quad (30)$$

Orthogonal condition is as

$$(\alpha+\beta-1)w_1 + (\beta-1)w_2 + (\beta-2)w_3 + (\beta-3)w_4 + (\alpha-1)w_5 - w_6 = 0 \quad (31)$$

and the primal-dual relations are as follows

$$AT^{\alpha+\beta-1} = w_1 d(w), B_1 T^{\beta-1} = w_2 d(w), CT^{\beta-2} = w_3 d(w)$$
$$, DT^{\beta-3} = w_4 d(w), ET^{\alpha-1} = w_5 d(w), FT^{-1} = w_6 d(w) \qquad (32)$$

Using the above primal-dual relation the followings are given by

$$\frac{B_1 w_1}{A w_2} = \left(\frac{C w_2}{B_1 w_3}\right)^{\alpha}, \left(\frac{B_1 w_1}{A w_2}\right)^{\frac{1}{\alpha}} = \left(\frac{D w_3}{C w_4}\right), \left(\frac{B_1 w_1}{A w_2}\right)^{\frac{(\beta-\alpha-2)}{\alpha}} = \frac{E w_4}{D w_5}, \frac{B_1 w_1}{A w_2} = \frac{F w_5}{E w_6} \qquad (33)$$

along with

$$T = \frac{C w_2}{B_1 w_3} \qquad (34)$$

Solving (30), (31) and (33) the critical value $w_1^*, w_2^*, w_3^*, w_4^*, w_5^*, w_6^*$ of the dual variables $w_1, w_2, w_3, w_4, w_5, w_6$ can be obtained and finally the optimum value $T^*$ of $T$ can be calculated from the equation of (34) substituting the critical values. Now the minimized total average cost $TOC^*_{\alpha,\beta}$ has been calculated by substituting $T^*$ in (28) analytically. The minimized total average cost and the optimal ordering interval is evaluated from (28) numerically.

**(ii) Case-2:** In this case we consider the rate of change of inventory level is fractional but integral memory index is absent *i.e.* $(\beta = 1.0, \ 0 < \alpha \leq 1.0)$.

In this case, the inventory model is

$$\begin{cases} \text{Min} TOC^{av}_{\alpha,1}(T) = AT^{\alpha} + B_1 T^0 + CT^{-1} + DT^{-2} + ET^{\alpha-1} \\ \text{Subject to } T \geq 0 \end{cases} \qquad (35)$$

$$A = \left(\frac{C_2 \gamma}{(\alpha+1)\Gamma(\alpha+1)}\right), B_1 = -\frac{C_2 \gamma t_1^{\alpha}}{\Gamma(\alpha+1)}, C = \frac{C_2 \gamma t_1^{\alpha+1}}{\Gamma(\alpha+1)} - \frac{C_2 \gamma t_1^{\alpha+1}}{(\alpha+1)\Gamma(\alpha+1)}, D = 0, E = \frac{P\gamma}{\Gamma(\alpha+1)}$$

$$F = \frac{C_1 \gamma t_1^{2\alpha+1}\left(B(\alpha+1,1) - B(2\alpha+1,1)\right)}{\Gamma(\alpha+1)} + C_3.$$

Using the similar analogy as previous, the minimized total average cost and optimal ordering interval has been solved from (35).

**(iii) Case-3:** Here, we consider presence of integral memory only *i.e.* rate of change of inventory level is unit but integral memory index present *i.e.* $(\alpha = 1.0, \ 0 < \beta \leq 1.0)$.

Then, the inventory model is as,

$$\begin{cases} \text{Min} TOC^{av}_{1,\beta}(T) = AT^{\beta} + B_1 T^{\beta-1} + CT^{\beta-2} + DT^{\beta-3} + ET^0 + FT^{-1} \\ \text{Subject to } T \geq 0 \end{cases} \qquad (36)$$

$$A = \left(\frac{C_2 \gamma}{2\Gamma(\beta)\Gamma(2)} - \frac{C_2(\beta-1)\gamma}{3\Gamma(\beta)\Gamma(2)}\right), B_1 = \left(\frac{C_2 \gamma(\beta-1)t_1}{2\Gamma(\beta)\Gamma(2)} - \frac{C_2 \gamma t_1}{\Gamma(\beta)\Gamma(2)}\right), C = \left(\frac{C_2 \gamma t_1^2}{\Gamma(\beta)\Gamma(2)} - \frac{C_2 \gamma t_1^2}{2\Gamma(\beta)\Gamma(2)}\right),$$

$$D = \frac{C_2 \gamma(\beta-1)t_1^3}{3\Gamma(\beta)\Gamma(2)} - \frac{C_2 \gamma(\beta-1)t_1^3}{2\Gamma(\beta)\Gamma(2)}, E = \frac{P\gamma}{\Gamma(2)}, F = \frac{C_1 \gamma\left(B(2,\beta) - B(3,\beta)\right)t_1^{2+\beta}}{\Gamma(\beta)\Gamma(2)} + C_3.$$

Using the similar way of case-1, the minimized total average cost and the optimal ordering interval has been solved from (36).

## 4. Numerical Example

**(a)** To illustrate numerically the developed fractional order inventory model, we consider empirical values of the various parameters in proper units as $P=20, C_3=250, C_1=2.5, C_2=1.2, \gamma=5, t_1=1.3456$ and required solution has been made using Matlab minimization method.

| $\alpha$ | $\beta$ | $T^*_{\alpha,\beta}$ | $TOC^*_{\alpha,\beta}$ |
|---|---|---|---|
| 0.1 | 1.0 | 167.3164 | 4.1316 |
| 0.2 | 1.0 | 81.5273 | 9.4674 |
| 0.3 | 1.0 | 51.6472 | 16.6348 |
| 0.4 | 1.0 | 35.8339 | 25.9524 |
| 0.5 | 1.0 | 25.6802 | 37.5886 |
| 0.6 | 1.0 | 18.3754 | 51.4457 |
| 0.7 | 1.0 | 12.7205 | 66.9596 |
| 0.8 | 1.0 | 8.1426 | 82.7853 |
| 0.9 | 1.0 | 4.4270 | 96.2934 |
| 1.0 | 1.0 | 1.8715 | 103.1555 |

**Table-3:** Optimal ordering interval and minimized total average cost for $\beta=1.0\ and\ 0<\alpha\leq 1.0$.

It is clear from the Table-3that in presence of differential memory index (i.e. rate of change of inventory level is affected by some exogenous effect such as environment of the society, quality of the product etc.)and in absence of integral memory index (i.e. the system affected by affected some exogenous effect like transportation system and shortage of the product those are not directly connected to the inventory level), the minimized total average cost is gradually decreases with increasing memory effect. In long memory effect, profit becomes high but business stay long time to reach the minimum value of the total average cost. Hence, moderate value of memory in the system is good for the business policy. Thus for good business policy, the planner should not carry too long previous memory and should consider the memory of a reasonable time period.

Now, when we fix the system in absence of differential memory index (see table-**4**) and allow the integral memory index then the minimized total average cost gradually decreases with gradually increasing memory effect for shortage started later (here for $t_1=1.3456$). But if shortage started quickly (for $t_1=0.3456$) then we see that the minimized total average cost becomes maximum at $\beta=0.5,$ then it gradually decreases below and above. In this case the optimal ordering interval is very low compare to the previous case. If the shortage start later then we see that minimized total average cost is high in long memory affected system and gradually decreases *i.e.* profit increases when integral memory increases. In this case the length of optimal ordering interval is long compare to the quickly shortage starting system but it not highly long as in the case of long differential memory system. moderate compared the short stock period.

Hence, from the concept of real market of inventory system and mathematical results conclude that for moderate i.e. (here long stock period) stock period is appropriate for the business.

| | | $t_1 = 1.3456$. | | $t_1 = 0.3456$. | |
|---|---|---|---|---|---|
| $\alpha$ | $\beta$ | $T^*_{\alpha,\beta}$ | $TOC^*_{\alpha,\beta}$ | $T^*_{\alpha,\beta}$ | $TOC^*_{\alpha,\beta}$ |
| 1.0 | 0.1 | 3.5850 | 100.5354 | 0.4144 | 100.1621 |
| 1.0 | 0.2 | 2.7919 | 101.0677 | 0.4062 | 100.2666 |
| 1.0 | 0.3 | 2.4492 | 101.5542 | 0.4001 | 100.3269 |
| 1.0 | 0.4 | 2.2540 | 101.9791 | 0.3956 | 100.3545 |
| 1.0 | 0.5 | 2.1278 | 102.3364 | 0.3921 | **100.3588** |
| 1.0 | 0.6 | 2.0403 | 102.6250 | 0.3895 | 100.3472 |
| 1.0 | 0.7 | 1.9772 | 102.8469 | 0.3875 | 100.3253 |
| 1.0 | 0.8 | 1.9307 | 103.0059 | 0.3861 | 100.2975 |
| 1.0 | 0.9 | 1.8964 | 103.1070 | 0.3853 | 100.2668 |
| 1.0 | 1.0 | 1.8715 | 103.1555 | 0.3848 | 100.2355 |

**Table-4:** Optimal ordering interval and minimized total average cost for $\alpha = 1.0$ and $0 < \beta \leq 1.0$.

In Table-5 we have presented the optimal ordering interval and minimized total average cost considering both memory effect simultaneously. It is clear from the Table-5 that the optimal ordering interval is long for this long memory affected system. The large optimal ordering interval implies that there may be some demurrage of inventory though the profit is high. Thus the memory affected model will be more realistic in the range $\alpha, \beta \in (0.6, 1.0)$. In reality if we consider more memory i.e. previous experiences then system will be disturbed due to its high restriction.

| (a) | | | | (b) | | | |
|---|---|---|---|---|---|---|---|
| $\alpha$ | $\beta$ | $T^*_{\alpha,\beta}$ | $TOC^*_{\alpha,\beta}$ | $\alpha$ | $\beta$ | $T^*_{\alpha,\beta}$ | $TOC^*_{\alpha,\beta}$ |
| 0.1 | 0.7 | $1.0000 \times 10^4$ | 0.4734 | 0.7 | 0.1 | $1.0000 \times 10^4$ | 7.0451 |
| 0.2 | 0.7 | $1.0000 \times 10^4$ | 1.6247 | 0.7 | 0.2 | $1.0000 \times 10^4$ | 7.4491 |
| 0.3 | 0.7 | $1.0000 \times 10^4$ | 4.4020 | 0.7 | 0.3 | $1.0000 \times 10^4$ | 8.8091 |
| 0.4 | 0.7 | 925.7379 | 9.6625 | 0.7 | 0.4 | $1.0000 \times 10^4$ | 12.9857 |
| 0.5 | 0.7 | 353.5743 | 18.0839 | 0.7 | 0.5 | $3.1723 \times 10^3$ | 24.1489 |
| 0.6 | 0.7 | 167.6143 | 30.7290 | 0.7 | 0.6 | 327.3122 | 37.2889 |
| 0.7 | 0.7 | 81.9894 | 48.2179 | 0.7 | 0.7 | 81.9894 | 48.2179 |
| 0.8 | 0.7 | 36.2569 | 69.8566 | 0.7 | 0.8 | 33.8214 | 56.5084 |
| 0.9 | 0.7 | 11.2350 | 91.6450 | 0.7 | 0.9 | 18.8792 | 62.5678 |
| 1.0 | 0.7 | 1.9772 | 102.8469 | 0.7 | 1.0 | 12.7205 | 66.9596 |

**Table-5-**Optimal ordering interval and minimized total average cost for
**(a)** $\beta = 0.7$ and $0 < \alpha \leq 1.0$. **(b)** $\alpha = 0.7$ and $0 < \beta \leq 1.0$.

In the next we will investigate the effect of differential memory index on positive inventory level with respect to time.

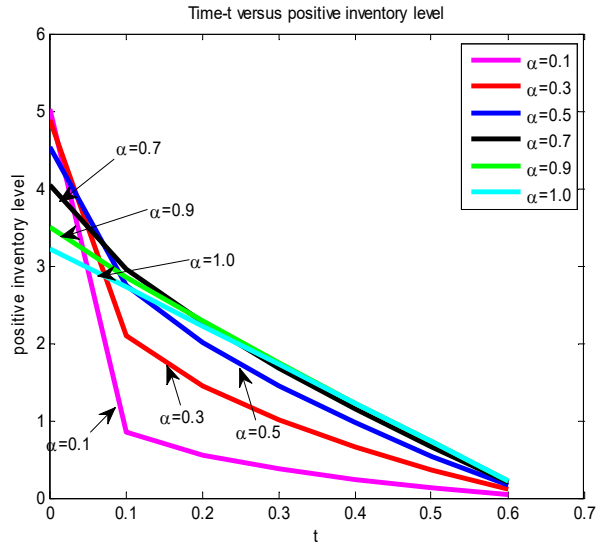

**Figure-2:** Positive inventory level versus time-*t* for different values of $\alpha$.

In figure-2 we have presented the positive inventory level with respect to *t* different values of the differential memory index ($\alpha$). It is clear from the figure that the positive inventory level changes linearly for low memory affected system but for long memory system then it initially fall rapidly and then decrease with slow rate

The graphical presentations of total average cost as a function of ordering interval (*T*) and shortage time ($t_1$) for different values of memory indices are presented in figure 3(a-p). It is clear figures that the nature of that total average cost is monotonic increasing function both the memory effect is high. But with the decrease of memory effect the average cost function become convex downwards. Hence the total minimized average cost will be minimized for intermediate value of the ordering interval and shortage time but in other cases the average cost function will be minimum at one end of the interval.

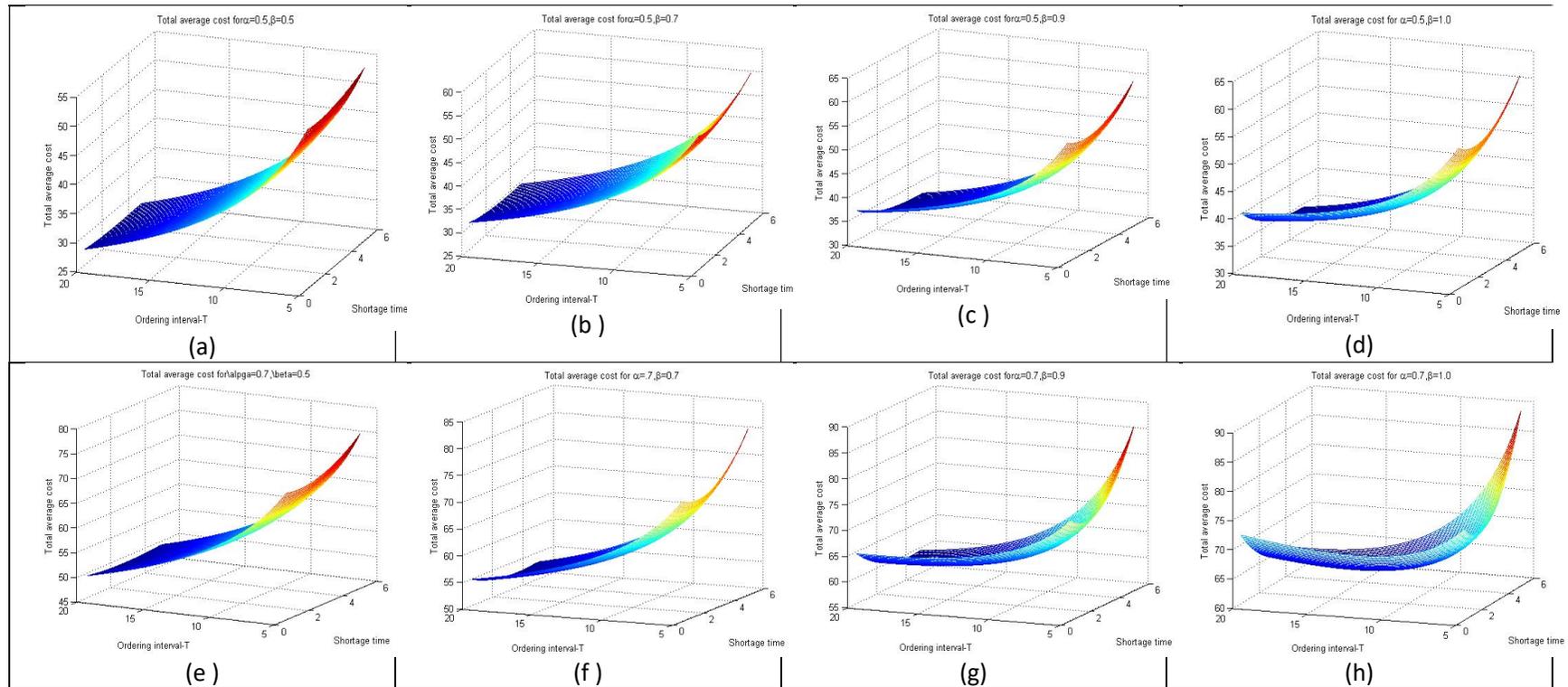

**Figure-3(a-h):** Total average cost versus ordering interval-*T* and shortage time-$t_1$ for different values of $\alpha, \beta$. (a-d) $\alpha = 0.5$ and $\beta = 0.5, 0.7, 0.9, 1.0$ (e-h) $\alpha = 0.7$ and $\beta = 0.5, 0.7, 0.9, 1.0$

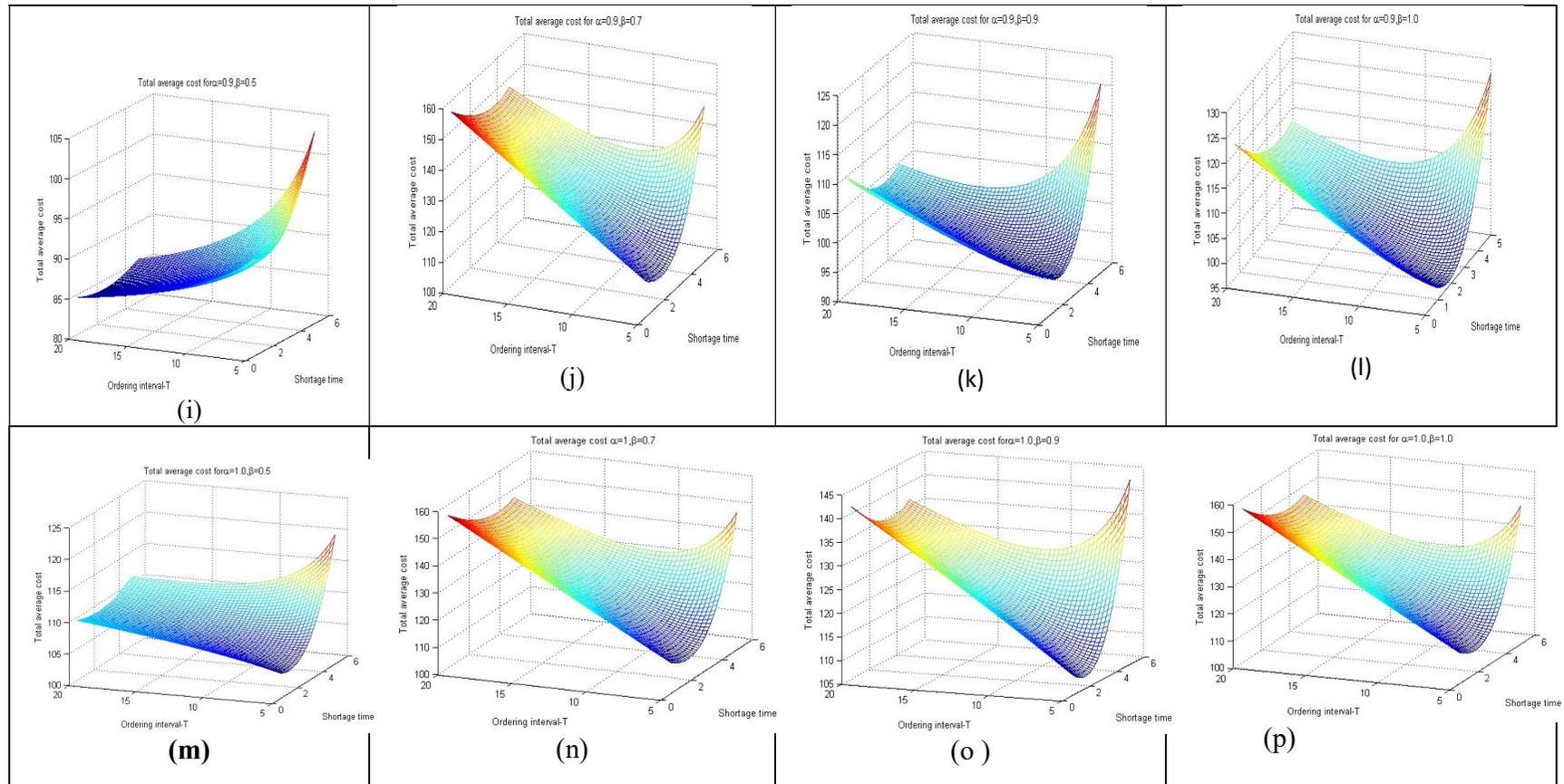

**Figure-3(i-p):** Total average cost versus ordering interval-$T$ and shortage time-$t_1$ for different values of $\alpha, \beta$. (i-l) $\alpha = 0.9$ and $\beta = 0.5, 0.7, 0.9, 1.0$ (m-p) $\alpha = 1$ and $\beta = 0.5, 0.7, 0.9, 1.0$

Next we will find the most important parameters for the considered system which has high influence in the total minimized average cost using sensitivity analysis

### 4.1 Sensitivity Analysis

The sensitivity analysis has been performed by changing each of the values of the considered parameter by +50%,+10%,-10%,-50% into show effects on the optimal ordering interval $T^*_{\alpha,1}$ and minimized total average cost $TOC^*_{\alpha,1}(T)$ taking one parameter $P, C_3, C_1, C_2, \gamma, t_1$ at a time and keeping the remaining parameters unchanged for the differential memory index $\alpha = 0.1$, and 0.9 *i.e.* in long memory effect or short memory effect.

| Parameter | Parameter Change (%) | $T^*_{\alpha,1}$ | $TOC^*_{\alpha,1}(T)$ | Parameter | Parameter Change(%) | $T^*_{\alpha,1}$ | $TOC^*_{\alpha,1}(T)$ |
|---|---|---|---|---|---|---|---|
| $P$ | +50% | 249.7255 | 4.5664 | $C_3$ | +50% | 167.3164 | 4.1316 |
| | +10%↑ | 183.7947 | 4.2320 | | +10%↑ | 167.3164 | 4.1316 |
| | -10% | 150.8405↑ | 4.0218↑ | | -10% | 167.3164 | 4.1316 |
| | -50% | 84.9788 | 3.4339 | | -50% | 167.3164 | 4.1316 |
| $\gamma$ | +50% | 167.3164 | 6.1974 | $t_1$ | +50% | 168.7113 | 3.8707 |
| | +10%↑ | 167.3164 | 4.5447↑ | | +10%↑ | 167.5878 | 4.0709↓ |
| | -10% | 167.3164 | 3.7184 | | -10% | 167.0494 | 4.1982 |
| | -50% | 167.3164 | 2.0658 | | -50% | ↑166.0350 | 4.5593 |
| $C_1$ | +50% | 168.0579 | 4.1358 | $C_2$ | +50% | 111.8959 | 5.5682 |
| | +10%↑ | 167.4647↑ | 4.1324 | | +10%↑ | 152.2021 | 4.4345 |
| | -10% | 167.1680 | 4.1307↑ | | -10% | 185.7889 | 3.8192↑ |
| | -50% | 166.5743 | 4.1273 | | -50% | 333.5486 | 2.4455 |

**Table-6:** The minimized total average cost and optimal ordering interval for $\alpha = 0.1$.

| Parameter | Parameter Change (%) | $T^*_{\alpha,1}$ | $TOC^*_{\alpha,1}(T)$ | Parameter | Parameter Change(%) | $T^*_{\alpha,1}$ | $TOC^*_{\alpha,1}(T)$ |
|---|---|---|---|---|---|---|---|
| $P$ | +50% | 5.9717 | 140.4055 | $C_3$ | +50% | 4.4270 | 96.2934 |
| | +10%↑ | 4.7268 | 105.2241 | | +10%↑ | 4.4270 | 96.2934 |
| | -10% | 4.1331 | 87.3026↑ | | -10% | 4.4270 | 96.2934 |
| | -50% | 3.0353↑ | 50.6686 | | -50% | 4.4270 | 96.2934 |
| $\gamma$ | +50% | 4.4270 | 144.4401 | $t_1$ | +50% | 5.4876 | 96.0662 |
| | +10%↑ | 4.4270 | 105.9227 | | +10%↑ | 4.6139 | 96.1332 |
| | -10% | 4.4270 | 86.6641↑ | | -10% | 4.2549 | 96.5165↓ |
| | -50% | 4.4270 | 48.1467 | | -50% | 3.7372 | 98.0877 |
| $C_1$ | +50% | 4.1603 | 96.8519 | $C_2$ | +50% | 3.3274 | 98.3645 |
| | +10%↑ | 4.4646 | 96.4069 | | +10%↑ | 4.1237 | 96.8128 |
| | -10% | 4.3890 | 96.1789↑ | | -10% | 4.8004 | 95.6972↑ |
| | -50% | 4.2315 | 95.7104 | | -50% | 7.8473 | 92.0014 |

**Table-7:** The minimized total average cost and optimal ordering interval for $\alpha = 0.9$.

It is clear from the Table-**6 and 7** that in long memory effect shortage cost per unit per unit time $(C_2)$ constant demand rate $(\gamma)$ are the most sensitive parameters but for low memory affected system the important parameters are per unit cost $(P)$ and constant demand rate $(\gamma)$.

## 5. Conclusion

In this paper, we have developed a fractional order EOQ model using the concept of memory dependency with the assumption that the demand is completely backlogged during the shortage time. To investigate the memory dependent EOQ model here we have introduced two type memory indexes as (i) differential memory index and (ii) integral memory index. In presence of differential memory index and in absence of integral memory index, profit becomes highest in long memory effect but optimal ordering interval is long. So, the presence of long memory is not practical in real life, the business man should consider moderate memory. It is clearly established from the numerical examples that in presence of integral memory index and in absence of differential memory index the profit gradually increases with increasing value of the memory effect. The sensitivity analysis shows that the parameter per unit cost of the total order quantity is the most sensitive parameter for the market studies in low memory effect but in long memory the most sensitive are parameters are shortage cost per unit and demand rate are the most sensitive parameter. In long memory effect, profit is always high compared to the low memory effect but the length of the ordering interval is long. Here, per unit cost is less sensitive in long memory affected system compare to low memory affected system. Hence, total order quantity needs to be adjusted for short past experience effect but not for long experience. For future research, the model can be extended for partial backlogging demand during shortage. This work will help the business policy maker to include the level of past experience in the marketing system.

**Acknowledgements**

The authors would like to thank the reviewers and the editor for the valuable comments and suggestions. The authors would also like to thank the Department of science and Technology, Government of India, New Delhi, for the financial assistance under AORC, Inspire fellowship Scheme towards this research work.